\documentclass[12pt]{amsart}
\usepackage{geometry}                
\usepackage{graphicx}
\usepackage{url}
\usepackage{amsmath}
\usepackage{amsthm}
\usepackage{amssymb}
\usepackage{epstopdf}
\usepackage{fancyhdr}

\usepackage{tikz} 
\usetikzlibrary{chains,shapes,arrows,trees,matrix,positioning,backgrounds,fit,calc,fadings}

\title{Reminiscences by a student of Langlands}
\author{Thomas Hales}

\begin{document}
\maketitle

\lhead{Hales}
\rhead{Reminiscences}

\parskip=0.8\baselineskip


\begin{quote}
  We are in a forest whose trees will not fall with a few
  timid hatchet blows.  We have to take up the double-bitted
  axe and the cross-cut saw, and hope that our muscles are
  equal to them.  -- R. P. Langlands
  \end{quote}

Bob Langlands was my thesis advisor at Princeton, 1984-1986.  No
mathematician has shaped my research career so profoundly as he has.
To use Weil's memorable phrase, I offer a few {\it souvenirs
  d'apprentissage}.

\section{Leap to Generality}

A stack exchange discussion asks for the largest
``leap-to-generality'' in mathematical history.  Suggestions include
the notion of category theory (Eilenberg Mac Lane), the rise of
abstract algebraic structures, Cantor's set theory, mathematics of the
infinite (starting with Archimedes' use of the method of exhaustion),
Aristotelian logic, the Turing machine, the foundations of probability
(Kolmogorov), and the Langlands program.

Harish-Chandra, Grothendieck, and Kolmogorov were Langlands's early
``models for emulation.''  In Langlands's own words, ``not satisfied
with partial insights and partial solutions, they [Harish-Chandra and
  Grothendieck] insisted -- not so much in the form of intentions or
exhortations as in what they brought to pass -- on methods that were
adequate to establishing the theories envisaged in their full natural
generality.''

\section{Princeton, 1983}

By the time I arrived in Princeton as a first-year graduate student in
the fall of 1983, I had already acquired interests including Lie theory,
representation theory, and the trace formula (thanks to Paul Cohen),
$p$-adic analysis (thanks to J.W.S. Cassels), and modular forms
(thanks to John Thompson in the heyday of moonshine).

From my first days at Princeton, I visited the Insititute for Advanced
Study (IAS) once or twice a week, attending the Borel-Mili\v ci\'c
$D$-modules seminar, Borel's 60th birthday conference, and the
Harish-Chandra memorial conference.

Eclipsing everything else that year were the momentous Morning and
Afternoon Seminars on the trace formula.  It is hard for me to convey
how deeply formative those seminars were for me, even if I was not
then at a stage to appreciate their full significance.  Experts --
especially Arthur, Clozel, and Rogawski -- encouraged me and taught me
the basics.  

I first met Bob Langlands in person in January 1984 at a dinner
arranged by Helaman Ferguson -- somebody that I spent considerable
time with my first year at Princeton. (Helaman's son, Samuel, later
became my coauthor on the proof of the Kepler conjecture.)  By that
spring, Langlands had become my advisor, and I had burrowed my way
into the Corvallis conference proceedings.  By arriving on scene after
the 1977 Corvallis conference, which was still the subject of spirited
conversations, I was made to feel I had missed a major event in the
history of the Langlands program.

My other main reference that first year was {\it Les D\'ebuts} (or the
purple turtle as we called it), where the the fundamental lemma was
first stated.  My research problem, broadly stated, was to use Igusa
theory to understand the transfer of $p$-adic orbital integrals
between a reductive group and its endoscopic groups.

This remained my primary research interest for ten years.  It has been
a great adventure to witness the trajectory of endoscopy over the
decades, culminating in Ng\^o Bao Ch\^au's proof of the fundamental
lemma.  In research posted to the arXiv last year, the fundamental
lemma has finally emerged in its natural geometric context, as
expressing that some dual abelian varieties have the same $p$-adic
volume.

\smallskip

Already by the time I met him, Langlands had a towering reputation for
his mathematical achievements.  He once published an unforgettable
critical book review with lede, ``This is a shallow book on deep
matters'' that compounded his formidable reputation.  I found that he
was more mellow in person than his reputation might suggest, and he
embodied the Institute's ideal of curiosity driven research.  It was
heartwarming for me to see Bob last year at the Abel Prize conference
in Minneapolis after many years.

\section{apprenticeship}

As a graduate student, the mathematical facts I learned mattered far
less than my apprenticeship as a researcher under Langlands.  I
arrived with good work habits, a disposition for long calculations,
and ambition.  Here are a few things my apprenticeship gave me.

\subsection{taste}
American popular culture failed miserably in conveying
great mathematical ideas to me.  In my teenage years, my (undeveloped)
idea of research mathematics was a confusing amalgamation of general
relativity, Thom's catastrophe theory, the Penrose staircase, and
stunning continued fraction expansions.  I remember wandering through
the library stacks and wondering which of these books matter most?

There is no question that my mathematical taste improved enormously
under Langlands. More broadly, lectures by Langlands, Serre, Weil,
Borel, Kottwitz, Iwasawa, Thurston, and Witten developed my tastes.

\subsection{seclusion}

Today, Google Scholar, the ArXiv, Wikipedia, and MathOverflow give us
nearly instant answers; polymath offers instant collaboration.

Then, there was a widespread belief that serious mathematical research
required long periods of intense work in relative seclusion.  In
Flexner's vision, the Institute ``exists as a paradise for scholars
who \ldots have won the right to do as they please and who accomplish
most when enabled to do so.''  I did not see Langlands summers, when
he went into work-related retreat in Montreal.  I developed my own
routines of seclusion: research retreats to a family cabin at the
Sundance ski resort in Utah, secret study areas, and long runs along
the Raritan Canal.

I have never seen Langlands at any mega-conference, and he did not
push the chores of professional service.  He discouraged rapid
publication.  He valued mathematical substance and frowned on veneer.
Details of proofs mattered.  He advised me to pick jobs based on
educational merit rather than salary or prestige.

\subsection{complexity}

There was a communal belief that we were building a monumental
edifice that would take many decades to complete.  One hundred page
research papers were the norm (and still are).

As documented in Wikipedia's list of long mathematical proofs, it is
no coincidence that some of the longest papers in mathematics are in
this field: Langlands (Eisenstein series), Arthur (trace formula),
Waldspurger (stable trace formula), Lafforgue (Langlands conjecture
for the general linear group over function fields); or in neighboring
fields in papers by Grothendieck, Hironaka, Harish-Chandra, Cartan, and Deligne.

As his student, I learned how to hold onto a problem that might take
years to solve.  I learned how to build evidence to support a hunch,
how to follow a lead, and how to bury a fruitless idea and move on.
These skills have transferred to my other large-scale research projects.

Researchers in the Langlands program were expected to rapidly
assimilate many research fields.  See Knapp's nine-page reading list
``Prerequisites for the Langlands Program,'' which expands to
thousands of pages of readings in algebraic geometry, Lie theory and
algebraic groups, representation theory, algebraic number theory, and
modular forms.

Yet book lists misguide us.  Library scholarship is not mathematical
research, and Langlands himself gave me just a few required
readings. When we met, it was entirely focused on what I was able to
calculate or figure out and on his suggestions for figuring things out
better.

\subsection{exoticism}

Langlands has described his recurrent dream of escaping -- in
T.E. Lawrence style -- ``into the life and language of some exotic
land and beginning anew,''
leading to his expeditions to Ankara.
When I was his student, he had just completed his
French lectures on stabilization of the trace formula and was
collaborating with Rapoport on their German masterpiece
``Shimuravariet\"aten und Gerben.''

Langlands's exoticism is inseparable from his mathematical oeuvre.  To
me, it explains how he went straight from a Yale thesis on PDEs and
analytic semigroups, to teaching a graduate course at Princeton on
class-field theory (``I still knew almost nothing about the subject,
had only two weeks to prepare, was very young, and scared stiff'').
With minimal ado, he would jump into an exotic field, then reconnect
it with the ever expanding Langlands program.  His new beginnings are
memorable, such as when he wrote that if his proofs seem clumsy, it
was because he ``has not cocycled before and has only minimum control
of his vehicle.''

\subsection{examples and generalization}

Langlands made many detailed studies of special cases, to shed light
on the general theory. There was base change for $GL(2)$,
Jacquet-Langlands for $GL(2)$, Labesse-Langlands for $SL(2)$, Igusa
theory and endoscopy for $SL(3)$, representations of abelian algebraic
groups, and a partial stable trace formula for $SL(n)$.

There were several axes of generalization.  What works for $SL(2)$ (or
even $GL(1)$) should work for all reductive groups.  What works for
the field of rational numbers should work for all global fields.  What
works for one local field should work for all.  There should be a
parallel between local and global theories, related by local-global
principles.

Methods were always under assessment: could they encompass the general
case?

\subsection{freedom}

Langlands's interests were already shifting to percolation theory when
I was his student. I was to be the last of his PhD students in the
Langlands program for decades.

The winter of the Langlands program lasted several years, starting
with Langlands's switch to percolation theory, and continuing until
Wiles's announcement of Fermat's Last Theorem.  Although significant
activity continued during the winter years, insiders and outsiders
alike had become increasingly disenchanted by the glacial pace towards
solutions to the central problems of the program.

At the time, Langlands would sometimes baffle his audience and say
speak on percolation theory to an audience clearly expecting
automorphic representations.  As his student, I have claimed the same
freedom to pursue my mathematical interests wherever they lead,
however baffling. On one level, I have left the Langlands program
behind.  But on another level, I have remained a true student of
Langlands by claiming this freedom.



\end{document}